\documentclass[12pt,reqno]{amsart}
\usepackage{enumerate, latexsym, amsmath, amsfonts, amssymb, amsthm, color}
\def\pmod #1{\ ({\rm{mod}}\ #1)}
\def\Z{\Bbb Z}

\def\t{\text}
\def\f{\frac}
 \def\Li{{\rm Li}}
\def\mo{{\rm{mod}\ }}

\def\ls{\leqslant}
\def\gs{\geqslant}

\def\bi{\binom}

\def\eq{\equiv}

\def\da{\delta}

\theoremstyle{plain}

\newtheorem{conjecture}{Conjecture}
\theoremstyle{definition}

\theoremstyle{remark}
\newtheorem{remark}{Remark}

 \vspace{4mm}

\begin{document}
 \baselineskip=16pt
\hbox{In: M. Kaneko, S. Kanemitsu and J. Liu (eds.), Number Theory: Plowing}
\hbox{and Starring through High Wave Forms, Proc. 7th China-Japan Seminar}
\hbox{(Fukuoka, Oct. 28--Nov. 1, 2013), Ser. Number Theory Appl., Vol. 11,}
\hbox{World Sci., Singapore, 2015, pp. 169--187.}
\medskip

\title
[Problems on combinatorial properties of primes]
{Problems on combinatorial properties of primes}

\author
[Zhi-Wei Sun] {Zhi-Wei Sun}

\thanks{Supported by the National Natural Science Foundation (grant 11171140)
 of China}

\address {Department of Mathematics, Nanjing
University, Nanjing 210093, People's Republic of China}
\email{zwsun@nju.edu.cn}

\keywords{Primes, combinatorial properties, conjectures, the prime-counting function, the $n$-th prime, partition functions.
\newline \indent 2010 {\it Mathematics Subject Classification}. Primary 11A41, 1B75; Secondary 05A10, 05A17, 11P32, 11P83.}

\begin{abstract} For $x\gs0$ let $\pi(x)$ be the number of primes not exceeding $x$. The asymptotic behaviors of the prime-counting function $\pi(x)$ and the $n$-th prime $p_n$
have been studied intensively in analytic number theory. Surprisingly, we find that $\pi(x)$ and $p_n$ have many combinatorial properties which should not be ignored.
In this paper we pose 60 open problems on combinatorial properties of primes (including connections between primes and partition functions)
for further research. For example,  we conjecture that for any integer $n>1$ one of the $n$ numbers
$\pi(n),\pi(2n),\ldots,\pi(n^2)$ is prime; we also conjecture that for any integer $n>6$ there exists a prime $p<n$
such that $pn$ is a primitive root modulo $p_n$. One of our conjectures involving the partition function $p(n)$ states that
for any prime $p$ there is a primitive root $g<p$ modulo $p$ with $g\in\{p(n):\ n=1,2,3,\ldots\}$.
\end{abstract}

\maketitle

\section{Introduction}
\setcounter{lemma}{0}
\setcounter{theorem}{0}
\setcounter{corollary}{0}
\setcounter{remark}{0}
\setcounter{equation}{0}
\setcounter{conjecture}{0}

  Prime numbers play important roles in number theory. For $x>0$ let $\pi(x)$ denote the number of primes not exceeding $x$. The celebrated Prime Number Theorem states that
 $$\pi(x)\sim\Li(x) \quad\t{as}\ x\to+\infty,$$
where $\Li(x)=\int_2^x\f{dt}{\log t}\sim \f x{\log x}$. This has the following equivalent version:
$$p_n\sim n\log n\quad\t{as}\ n\to+\infty,$$
where $p_n$ denotes the $n$-th prime. To get sharp estimations for $\pi(x)$ is a main research topic in analytic number theory.
It is known (cf. \cite{Sc}) that under Riemann's Hypothesis we have
$$\pi(x)=\Li(x)+O(\sqrt x\log x)\quad\t{and}\quad p_{n+1}-p_n=O(\sqrt{p_n}\log p_n).$$
For convenience, we also set $\pi(0)=0$.

Many number theorists generally consider primes irregular and only focus on their asymptotic behaviours.
In contrast with the great achievements on the asymptotic behaviors of $\pi(x)$ and $p_n$ (see \cite{Z} for a recent breakthrough on prime gaps),
almost nobody has investigated combinatorial properties of primes
seriously and systematically.

Surprisingly, we find that the functions $\pi(x)$ and $p_n$ have many unexpected combinatorial properties depending on their exact values.
Also, partition functions arising from combinatorics
have nice connections with primes.
In this paper we pose 60 typical conjectures in this direction.
The next section contains 25 conjectures on combinatorial properties of $\pi(x)$, while Section 3 contains 25 conjectures on combinatorial properties of the function $p_n$.
Section 4 is devoted to 10 conjectures on primes related to partition functions. The reader may also consult \cite{S13a}
for the author's previous conjectures on alternating sums of consecutive primes.

The 60 selected conjectures in Sections 2--4 are somewhat incredible. Nevertheless, our numerical computations and related graphs in \cite{S} provide strong evidences to support them.
The author would like to offer 1000 Chinese dollars as the prize for the first complete solution to any one of the 60 conjectures.
We hope that the problems here might interest some number theorists and stimulate further research, but the solutions to most of them might be beyond the intelligence of human beings.

\section{Combinatorial properties of $\pi(x)$ and related things}
\setcounter{lemma}{0}
\setcounter{theorem}{0}
\setcounter{corollary}{0}
\setcounter{remark}{0}
\setcounter{equation}{0}
\setcounter{conjecture}{0}

\begin{conjecture}\label{Conj2.1} {\rm (2014-02-09)} {\rm (i)} For any integer $n>1$, $\pi(kn)$ is prime for some $k=1,\ldots,n$.
Moreover, for every $n=1,2,3,\ldots$, there is a positive integer $k<3\sqrt n+3$ with $\pi(kn)$ prime.

{\rm (ii)} Let $n_0=5,\ n_1=3$ and $n_{-1}=6$. For each $\da\in\{0,\pm1\}$ and any integer $n>n_{\da}$, there is a positive integer $k<n$ such that $k^2+k-1$ and $\pi(kn)+\da$ are both prime.
\end{conjecture}

\begin{remark}\label{Rem2.1} (a) We also conjecture that for any integer $n>92$ there is a prime $p\ls n$ with $\pi(pn)$ prime.
We have verified part (i) of Conjecture \ref{Conj2.1} for $n$ up to $2\times10^7$, and our data and graphs for the sequence $a(n)=|\{0<k<n: \ \pi(kn)\ \t{is prime}\}|\ (n=1,2,3,\ldots)$
(cf. \cite[A237578]{S})
strongly support its truth. It seems that $|\{1\ls k\ls n:\ \pi(kn)\ \t{is prime}\}|\sim \pi(n)/2$ as $n\to+\infty$.
See also \cite[A237615]{S} for part (ii) of Conjecture \ref{Conj2.1},
and note that it is not yet proven that there are infinitely many primes of the form $x^2+x-1$ with $x\in\Z$.

(b) We also conjecture the following analogue of Conjecture \ref{Conj2.1} (cf. \cite[A238703]{S}):
For any integers $n>m>0$ with $m\nmid n$, there is a positive
integer $k<n$ with $\lfloor kn/m\rfloor$ prime. Note that $\lfloor kn/m\rfloor$ is the number of multiples of $m$ among $1,2,\ldots,kn$.
\end{remark}

\begin{conjecture}\label{Conj2.2} {\rm (2014-02-10)} {\rm (i)} For any positive integer $n$, there is a positive integer $k<p_n$ such that $\pi(kn)\eq0\pmod n$.

{\rm (ii)} For each positive integer $n$, the set $\{\pi(kn):\ k=1,\ldots,2p_n\}$ contains a complete system of residues modulo $n$.
\end{conjecture}

\begin{remark}\label{Rem2.2} See \cite[A237597]{S} and \cite[A237643]{S} for related sequences concerning this conjecture.
\end{remark}

\begin{conjecture}\label{Conj2.3} {\rm (2014-02-20)} Let $n>1$ be an integer. Then $\pi(jn)\mid \pi(kn)$ for some $1\ls j<k\ls n$ with $k\eq1\pmod j$.
\end{conjecture}

\begin{remark}\label{Rem2.3} For example, $\pi(3\times 50) =35$ divides $\pi(7\times50)=70$ with $7\eq1\ (\mo\ 3)$.
We have verified the conjecture for all $n=2,3,\ldots,30000$. See \cite[A238224]{S} for a related sequence.
\end{remark}

\begin{conjecture}\label{Conj2.4} {\rm (i) (2014-02-10)} For any positive integer $n$, there is a positive integer $k<p_n$ such that $\pi(kn)$ is a square.

{\rm (ii) (2014-02-14)} Let $n$ be any positive integer. Then, for some $k=1,\ldots,n$, the number of twin prime pairs not exceeding $kn$ is a square.
\end{conjecture}

\begin{remark}\label{Rem2.4} See \cite[A237598, A237612, A237840, A237879 and A237975]{S} for some sequences related to this conjecture.
Similar to part (i), we conjecture that for any integer $n>9$ there is a positive integer $k<p_n/2$ such that
$\pi(kn)$ is a triangular number. We have verified part (ii) of the conjecture for all $n=1,\ldots,22000$; for example, for $n=19939$ we may take $k=12660$ since there are exactly
$1000^2=10^6$ twin prime pairs not exceeding $12660\times 19939=252427740$.
\end{remark}

\begin{conjecture}\label{Conj2.5}  {\rm (i) (2014-02-24)} For any integer $n>5$, there is a positive integer $k<n$ with $kn+\pi(kn)$ prime.

{\rm (ii) (2014-03-06)} If $n$ is a positive integer, then $p_{kn}-\pi(kn)$ is prime for some $k=1,\ldots,n$.
\end{conjecture}

\begin{remark}\label{Rem2.5} See \cite[A237712 and A238890]{S} for related sequences. Part (ii) of the conjecture implies that there are infinitely many primes $p$
with $p-\pi(\pi(p))$ prime.
\end{remark}

\begin{conjecture}\label{Conj2.6} {\rm (2014-03-07)} {\rm (i)} For any integer $n>2$, there is a prime $p\ls n$ with $\pi(\pi((p-1)n))$ prime.

{\rm (ii)} Let $n$ be any positive integer. Then $\pi(\pi(kn))$ is a square for some $k=1,\ldots,n$. Also, there exists a positive integer $k\ls(n+1)/2$
such that $\pi(\pi(kn))$ is a triangular number.
\end{conjecture}

\begin{remark}\label{Rem2.6} See \cite[A238504, A238902 and A239884]{S} for related data and graphs. We have verified the two assertions in Conjecture \ref{Conj2.6}(ii) for $n$ up to $2\times10^5$ and $10^5$ respectively; for example,
\begin{align*}\pi(\pi(8514\times9143))&=\pi(4550901)=565^2,
\\\pi(\pi(37308\times98213))&=\pi(174740922)=3123^2,
\\\pi(\pi(83187\times192969))&=\pi(715034817) = 6082^2.
\end{align*}
We guess that there are positive constants $c_1$ and $c_2$ such that
$$c_1\ls\f{|\{1\ls k\ls n:\ \pi(\pi(kn))\ \t{is a square}\}|}{\log n}\ls c_2\qquad\t{for all}\ n=2,3,\ldots.$$
\end{remark}

\begin{conjecture}\label{Conj2.7}  {\rm (i) (2014-02-17)} For any integer $n>4$ and $k=1,\ldots,n$, we have $\pi(kn)^{1/k}>\pi((k+1)n)^{1/(k+1)}$.

{\rm (ii) (2014-02-22)} Let $n$ be any positive integer. Then, $\pi((k+1)n)-\pi(kn)\ ($the number of primes in the interval $(kn,(k+1)n])$ is a square for some $k=0,\ldots,n-1$.
\end{conjecture}

\begin{remark}\label{Rem2.7} For any integer $n>1$, Bertrand's postulate (first proved by Chebyshev in 1850) indicates that $\pi(2n)>\pi(n)$, and
Oppermann's conjecture states that $\pi((n-1)n)<\pi(n^2)<\pi(n(n+1))$. Our computation suggests that $\pi(kn)<\pi((k+1)n)$ for any integers $n\gs k>0$.
 A conjecture of Firoozbakht (cf. \cite[p.\,185]{R}) asserts that the sequence $p_n^{1/n}\ (n=1,2,3,\ldots)$ is strictly decreasing,
 and the author's recent paper \cite{S13b} contains many similar conjectures on monotonicity of arithmetical sequences.
See also \cite[A238277]{S} for a sequence related to Conjecture \ref{Conj2.7}(ii).
\end{remark}

\begin{conjecture}\label{Conj2.8} {\rm (i) (2014-03-16)} Let $n>3$ be an integer. Then $\pi(pn)-\pi((p-1)n)$ is prime for some prime $p<n$.
Also, there is an odd prime $p\ls n$ with $\pi(\f{p+1}2n)-\pi(\f{p-1}2n)$ prime.

{\rm (ii) (2014-02-22)} For any integer $n>3$, there is a number $k\in\{1,\ldots,n-1\}$ with $\pi(kn)-\pi((k-1)n)$ and $\pi((k+1)n)-\pi(kn)$ both prime.
\end{conjecture}

\begin{remark}\label{Rem2.8} See \cite[A239328, A239330 and A238278]{S} for related data and graphs.
\end{remark}

\begin{conjecture}\label{Conj2.9} {\rm (2014-02-22)}
{\rm (i)} For any integer $n>1$, there is a positive integer $k<n$ such that the intervals $(kn,(k+1)n)$ and $((k+1)n,(k+2)n)$ contain the same number of primes, i.e.,
 $$\pi(kn),\ \pi((k+1)n),\ \pi((k+2)n)$$ form a three-term arithmetic progression.

{\rm (ii)} For any integer $n>4$, there is a positive integer $k<p_n$ such that
$$\pi(kn),\ \pi((k+1)n),\ \pi((k+2)n),\ \pi((k+3)n)$$ form a four-term arithmetic progression.
\end{conjecture}

\begin{remark}\label{Rem2.9} See \cite[A238281]{S} for a sequence related to part (i) of this conjecture.
\end{remark}

\begin{conjecture}\label{Conj2.10} {\rm (2014-02-23)} For any positive integer $n$, we have
$$|\{\pi((k+1)n)-\pi(kn):\ k=0,\ldots,n-1\}|\gs\sqrt{n-1},$$
and equality holds only when $n$ is $2$ or $26$.
\end{conjecture}

\begin{remark}\label{Rem2.10} See \cite[A230022]{S} for related data and graphs.
\end{remark}

\begin{conjecture}\label{Conj2.11} {\rm (2014-02-24)} Let $n>1$ be an integer. Then, for some prime $p\ls p_n$, the three numbers $\pi(p),\pi(p+n),\pi(p+2n)$ form a nontrivial arithmetic progression,
i.e., $\pi(p+2n)-\pi(p+n)=\pi(p+n)-\pi(p)>0$.
\end{conjecture}

\begin{remark}\label{Rem2.11} See \cite[A210210]{S} for a related sequence.
\end{remark}

\begin{conjecture}\label{Conj2.12} {\rm (i) (2014-02-08)} Each integer $n>10$ can be written as $k+m$ with $k$ and $m$ positive integers such that $\pi(km)\ (\t{or}\ \pi(k^2m))$ is prime.

{\rm (ii) (2014-03-20)} For any integer $n>9$, there are positive integers $k$ and $m$ with $k+m=n$ such that $\pi(2k)-\pi(k)$ and $\pi(2m)-\pi(m)$ are both prime. Also,
any integer $n>4$ can be written as a sum of two positive integers $k$ and $m$ such that $\pi(2k)-\pi(k)$ is a prime and $\pi(2m)-\pi(m)$ is a square.
\end{conjecture}

\begin{remark}\label{Rem2.12} See \cite[A237497, A239428 and A239430]{S} for related data and graphs.
As $\pi(2k+2)-\pi(k+1)-(\pi(2k)-\pi(k))\in\{0,\pm1\}$ and $\pi(2n)-\pi(n)\sim n/\log n$, there are infinitely many positive integers $k$ with $\pi(2k)-\pi(k)$ prime.
Similar to part (i), we also conjecture (cf. \cite[A237531]{S}) that for any integer $n>5$
there is a positive integer $k<n/2$ such that $\varphi(k(n-k))-1$ and $\varphi(k(n-k))+1$ are twin prime, where $\varphi$ denotes Euler's totient function.
\end{remark}

For $x\gs0$, we use $\pi_2(x)$ to denote the number of twin prime pairs not exceeding $x$, i.e., $\pi_2(x)=|\{p\ls x:\ p\ \t{and}\ p-2\ \t{are both prime}\}|$.

\begin{conjecture}\label{Conj2.13} {\rm (2014-02-15)} {\rm (i)} Each integer $n>5$ can be written as $k+m$ with $k$ and $m$ positive integers such that $\pi_2(km)$ is prime.

{\rm (ii)} Any integer $n>8$ can be written as $k+m$ with $k$ and $m$ positive integers such that $\pi_2(km)-1$ and $\pi_2(km)+1$ are twin prime.
\end{conjecture}

\begin{remark}\label{Rem2.13} This is an analogue of Conjecture \ref{Conj2.12}(i) for twin prime pairs.
\end{remark}

\begin{conjecture}\label{Conj2.14} {\rm (i) (2014-02-11)} For any positive integer $n$, the set $\{\pi(k^2):\ k=1,\ldots,2p_{n+1}-3\}$ contains a complete system of residues modulo $n$.

{\rm (ii) (2014-02-17)} The sequence $\root n\of{\pi(n^2)}\ (n=3,4,\ldots)$ is strictly decreasing.

{\rm (iii) (2014-02-17)} For any integer $n>0$, the interval $[\pi(n^2),\pi((n+1)^2)]$ contains at least one prime except for $n=25,\,35,\,44,\,46,\,105$.
\end{conjecture}

\begin{remark}\label{Rem2.14} Legendre's conjecture asserts that for each positive integer $n$ there is a prime between $n^2$ and $(n+1)^2$.
\end{remark}

\begin{conjecture}\label{Conj2.15} {\rm (i) (2014-02-11)} For any integer $n>8$, $\pi(k)$ and $\pi(k^2)$ are both prime for some integer $k\in(n,2n)$.

{\rm (ii) (2014-02-11)} There are infinitely many primes $p$ with $\pi(p)$, $\pi(\pi(p))$ and $\pi(p^2)$ all prime.

{\rm (iii) (2014-04-09)} For any positive integer $n$, there are infinitely many primes $p$ with $\pi(kp)$ prime for all $k=1,\ldots,n$.
\end{conjecture}

\begin{remark}\label{Rem2.15} See \cite[A237657, A237687 and A240604]{S} for related sequences and data.
\end{remark}

\begin{conjecture}\label{Conj2.16} {\rm (i) (2014-02-09)} For any integer $n>1$,  $\pi(n+k^2)$ is prime for some $k=1,\ldots,n-1$. In general, for each $a=2,3,\ldots$, if an integer $n$ is sufficiently
large, then $\pi(n+k^a)$ is prime for some $k=1,\ldots,n-1$.

{\rm (ii) (2014-02-10)} Let $n>4$ be an integer. Then $n+\pi(k^2)$ is prime for some $k=1,\ldots,n$.

{\rm (iii) (2014-03-01)} If a positive integer $n$ is not a divisor of $12$, then $n^2+\pi(k^2)$ is prime for some $1<k<n$.
For any integer $n>4$, $\pi(n^2)+\pi(k^2)$ is prime for some $1<k<n$. Also, for each $n=2,3,\ldots$ there is a positive integer $k<n$ such that
$\pi((k+1)^2)-\pi(k^2)$ and $\pi(n^2)-\pi(k^2)$ are both prime.
\end{conjecture}

\begin{remark}\label{Rem2.16}  See \cite[A237582, A237595 and A238570]{S} for related sequences. In 2012 the author conjectured that if $n$ is a positive integer then $n+k$ and $n+k^2$ are both prime for some $k=0,\ldots,n$ (cf. \cite[A185636]{S}).
\end{remark}

\begin{conjecture}\label{Conj2.17} {\rm (i) (2014-02-08)} For any integer $n>4$, there is a prime $p<n$ with $pn+\pi(p)$ prime. Moreover, for every positive integer $n$,
there is a prime $p<\sqrt{2n}\log(5n)$ with $pn+\pi(p)$ prime.

{\rm (ii) (2014-03-02)} For any integer $n>2$, there is a prime $p\ls n$ with $2\pi(p)-(-1)^n$ and $pn+((-1)^n-3)/2$ both prime.
\end{conjecture}

\begin{remark}\label{Rem2.17} See \cite[A237453 and A238643]{S} for related data and graphs. We have verified parts (i) and (ii) for $n$ up to $10^8$.
As a supplement to part (i), we also conjecture that for every $n=1,2,3,\ldots$ there is a positive integer $k<3\sqrt n$ such that $kn+p_k=\pi(p_k)n+p_k$ is prime.
Part (ii) implies that for any odd prime $p$ there is a prime $q\ls p$ with $pq-2$ prime.  By Chen's work \cite{C}, there are infinitely many
primes $p$ with $p+2$ a product of at most two primes.
\end{remark}

\begin{conjecture}\label{Conj2.18} {\rm (2013-11-24)} {\rm (i)}  Every $n=4,5,\ldots$ can be written as $p+q-\pi(q)$, where $p$ and $q$ are odd primes not exceeding $n$.

{\rm (ii)} For any integer $n>7$, there is a prime $p<n$ with $n+p-\pi(p)$ also prime.
\end{conjecture}

\begin{remark}\label{Rem2.18} We have verified part (i) for all $n=4,5,\ldots,10^8$; for example, $9=7+5-\pi(5)$ with $7$ and $5$ prime.
See \cite[A232463 and A232443]{S} for related sequences.
\end{remark}

\begin{conjecture}\label{Conj2.19} {\rm (2014-02-06)} {\rm (i)} For any integer $n>2$, there is a prime $p<2n$ with $\pi(p)$ and $2n-p$ both prime.

{\rm (ii)} For any integer $n>36$, we can write $2n-1=a+b+c$ with $a,b,c$ in the set $\{p:\ p\ \t{and}\ \pi(p)\ \t{are both prime}\}$.
\end{conjecture}

\begin{remark}\label{Rem2.19} Part (i) is a refinement of Goldbach's conjecture, and part (ii) is stronger than the weak Goldbach conjecture finally proved by Helfgott \cite{He}.
See \cite[A237284 and A237291]{S} for related representation functions.
\end{remark}

Recall that a prime $p$ with $2p+1$ also prime is called a Sophie Germain prime.

\begin{conjecture}\label{Conj2.20} {\rm (2014-02-13)} {\rm (i)} For any integer $n>4$, there is a prime $p<n$ such that $\pi(n-p)$ is a Sophie Germain prime.
Also, for any integer $n>8$ there is a prime $p<n$ such that $\pi(n-p)-1$ and $\pi(n-p)+1$ are twin prime.

{\rm (ii)} For any integer $n>4$, there is a prime $p<n$ such that $3m\pm1$ and $3m+5$ are all prime with $m=\pi(n-p)$. Also, for any integer $n>8$,
there is a prime $p<n$ such that $3m\pm1$ and $3m-5$ are all prime with $m=\pi(n-p)$.
\end{conjecture}

\begin{remark}\label{Rem2.20} See \cite[A237768 and A237769]{S} for related sequences. We have verified part (i) for $n$ up to $2\times 10^7$.
\end{remark}

\begin{conjecture}\label{Conj2.21} {\rm (2014-02-13)} {\rm (i)} For any integer $n>4$, there is a prime $p<n$ such that the number of Sophie Germain primes among $1,\ldots,n-p$
is a Sophie Germain prime.

{\rm (ii)} For any integer $n>12$, there is a prime $p<n$ such that
$$r=|\{q\ls n-p:\ q\ \t{and}\ q+2\ \t{are twin prime}\}|$$
and $r+2$ are twin prime.
\end{conjecture}

\begin{remark}\label{Rem2.21} See \cite[A237815 and A237817]{S} for related sequences.
\end{remark}

\begin{conjecture}\label{Conj2.22} {\rm (2014-02-12)} {\rm (i)} For any integer $n>2$, there is a prime $p<n$ such that $\pi(n-p)$ is a square.
Also, for any integer $n>2$ there is a prime $p<n$ such that $\pi(n-p)$ is a triangular number.

{\rm (ii)} For any integer $n>2$, there is a prime $p\ls p_n$ such that $\pi(n+p)$ is a square.
\end{conjecture}

\begin{remark}\label{Rem2.22} See \cite[A237706 and A237710]{S} for related sequences.
We have verified the first assertion in part (i) for $n$ up to $5\times10^8$, and guessed that the number of primes $p<n$ with $\pi(n-p)$ a square
is asymptotically equivalent to $c\sqrt n$ with $c$ a constant in the interval $(0.2,0.22)$. We also conjecture the following analogue (cf. \cite[A238732]{S})
of Conjecture \ref{Conj2.22}(i):
For any integers $m>2$ and $n>2$, there is a prime $p<n$ such that $\lfloor (n-p)/m\rfloor$ is a square.
\end{remark}

\begin{conjecture}\label{Conj2.23} {\rm (2014-02-13)} {\rm (i)} For any integer $n>11$, there is a prime $p<n$ such that the number of Sophie Germain primes among $1,\ldots,n-p$
is a square.

{\rm (ii)} For any integer $n\gs54$, there is a prime $p<n$ such that the number of Sophie Germain primes among $1,\ldots,n-p$ is a cube.
\end{conjecture}

\begin{remark}\label{Rem2.23} Part (i) is an analogue of Conjecture 2.22(i) for Sophie Germain primes. See \cite[A237837]{S} for a sequence related to part (ii).
\end{remark}

\begin{conjecture}\label{Conj2.24} {\rm (2014-03-02)} {\rm (i)} For every $n=2,3,\ldots$, there is an odd prime $p<2n$ such that the number of squarefree integers among $1,\ldots,\f{p-1}2n$
is prime.

{\rm (ii)} For any integer $n>3$, there is a prime $p<n$ such that the number of squarefree numbers among $1,\ldots,n-p$ is prime.
\end{conjecture}

\begin{remark}\label{Rem2.24}  See \cite[A238645 and A238646]{S} for related sequences.
\end{remark}

\begin{conjecture}\label{Conj2.25} {\rm (2014-02-22)} For any integer $n>4$, there is a number $k\in\{1,\ldots,n\}$ such that the number of prime ideals of
the Gaussian ring $\Z[i]$ with norm not exceeding $kn$ is a prime congruent to $1$ modulo $4$.
\end{conjecture}

\begin{remark}\label{Rem2.25} $\Z[i]$ is a principal ideal domain, and any prime ideal $P$ of it has the form $(p)$ with $p$ a rational prime congruent to $3$ modulo $4$ or
$p=a+bi$ with $N(p)=a^2+b^2$ a rational prime not congruent to $3$ modulo $4$. (Cf. \cite[p.\,120]{IR}.) So, the number of prime ideals of $\Z[i]$ with norm not exceeding $x$ actually equals
$$\pi(\sqrt x)+|\{\sqrt x<p\ls x:\ p\ \t{is a prime with}\ p\not\equiv 3\ (\mo\ 4)\}|.$$
\end{remark}

\section{Combinatorial properties involving the function $p_n$}
\setcounter{lemma}{0}
\setcounter{theorem}{0}
\setcounter{corollary}{0}
\setcounter{remark}{0}
\setcounter{equation}{0}
\setcounter{conjecture}{0}

\begin{conjecture}\label{Conj3.1} {\rm (Unification of Goldbach's Conjecture and the Twin Prime Conjecture, 2014-01-29)} For any integer $n>2$,
there is a prime $q$ with $2n-q$ and $p_{q+2}+2$ both prime.
\end{conjecture}

\begin{remark}\label{Rem3.1} We have verified this for $n$ up to $2\times 10^8$. See \cite[A236566]{S} for a related sequence.
Note that the conjecture implies the Twin Prime Conjecture. In fact, if all primes $q$ with $p_{q+2}+2$ prime are smaller than an even number $N>2$, then
for any such a prime $q$ the number $N!-q$ is composite since $N!-q\eq0\pmod{q}$ and $N!-q\gs q(q+1)-q>q$.
\end{remark}

\begin{conjecture}\label{Conj3.2} {\rm (Super Twin Prime Conjecture, 2014-02-05)} Any integer $n>2$ can be written as $k+m$ with $k$ and $m$ positive integers such that
$p_k+2$ and $p_{p_m}+2$ are both prime.
\end{conjecture}

\begin{remark}\label{Rem3.2} We have verified the conjecture for $n$ up to $10^9$. See \cite[A218829, A237259, A237260]{S} for related sequences.
If $p, p+2$ and $\pi(p)$ are all prime, then we call $\{p,p+2\}$ a {\it super twin prime} pair. Conjecture \ref{Conj3.2} implies that there are infinitely many
super twin prime pairs. In fact, if all those positive integers $m$ with $p_{p_m}+2$ prime are smaller than an integer $N>2$, then by Conjecture \ref{Conj3.2}, for each $j=1,2,3,\ldots$,
there are positive integers $k(j)$ and $m(j)$ with $k(j)+m(j)=jN$ such that $p_{k(j)}+2$ and $p_{p_{m(j)}}+2$ are both prime, and hence $k(j)\in((j-1)N,jN)$ since $m(j)<N$; thus
$$\sum_{j=1}^{\infty}\f1{p_{k(j)}}\gs\sum_{j=1}^\infty\f1{p_{jN}},$$
which is impossible since the series on the right-hand side diverges while the series on the left-hand side converges by Brun's theorem on twin primes (cf. \cite[p.\, 14]{CP}).
\end{remark}

\begin{conjecture}\label{Conj3.3} {\rm (2014-01-28)} Any integer $n>2$ can be written as $k+m$ with $k$ and $m$ positive integers such that
both $\{6k\pm1\}$ and $\{p_m,p_m+2\}$ are twin prime pairs.
\end{conjecture}

\begin{remark}\label{Rem3.3} Clearly this implies the Twin Prime Conjecture. We have verified Conjecture 3.3 for $n$ up to $2\times 10^7$. See \cite[A236531]{S} for related data and graphs.
\end{remark}

\begin{conjecture}\label{Conj3.4} {\rm (2014-02-07)} Any integer $n>1$ can be written as $k+m$ with $k$ and $m$ positive integers such that
$p_k^2-2,\ p_m^2-2$ and $p_{p_m}^2-2$ are all prime.
\end{conjecture}

\begin{remark}\label{Rem3.4} We have verified this for $n$ up to $10^8$. See \cite[A237413 and A237414]{S} for related sequences. It is not yet proven that there are infinitely many primes of the form $x^2-2$ with $x\in\Z$.
\end{remark}

\begin{conjecture}\label{Conj3.5} {\rm (i) (2013-11-25)} The set
$$\{k+m:\ 0 <k<m<n,\ \t{and}\ p_k, p_m, p_n \ \t{form an arithmetic progression}\}$$
coincides with $\{5,6,7,\ldots\}$.

{\rm (ii) (2014-02-22)} For every $n=1,2,3,\ldots$, there is a positive integer $k\ls 3p_n+8$ such that
$p_{kn},p_{(k+1)n},p_{(k+2)n}$ form a three-term arithmetic progression.
\end{conjecture}

\begin{remark}\label{Rem3.5} Recall that the Green-Tao theorem (cf. \cite{GT}) asserts that there are arbitrarily long arithmetic progressions of primes.
See \cite[A232502 and A238289]{S} for related sequences.
\end{remark}

\begin{conjecture}\label{Conj3.6} {\rm (2014-03-01)} {\rm (i)} For any integer $n>6$, there is a number $k\in\{1,\ldots,n\}$ with $p_{kn}+2\ (\t{or}\ p_{k^2n}+2)$ prime.
Moreover, for every $n=1,2,3,\ldots$ there is a positive integer $k<3\sqrt n+6$ with $p_{kn}+2$ prime.

{\rm (ii)} For any positive integer $n$, there is a number $k\in\{1,\ldots,n\}$ such that $2k+1$ and $p_{kn}^2-2$ are both prime.
\end{conjecture}

\begin{remark}\label{Rem3.6} Clearly part (i) is stronger than the Twin Prime Conjecture, while part (ii) implies that there are infinitely many primes $p$ with $p^2-2$ prime.
See \cite[A238573 and A238576]{S} for related data and graphs.
\end{remark}

\begin{conjecture}\label{Conj3.7} {\rm (i) (2013-12-01)} There are infinitely many positive integers $n$ such that
$$n\pm1,\ p_n\pm n,\ np_n\pm1$$
are all prime.

{\rm (ii) (2014-01-20)} There are infinitely many primes $q$ with $p_q^2+4q^2$ and $q^2+4p_q^2$ both prime.
\end{conjecture}

\begin{remark}\label{Rem3.7} For part (i), the first such a number $n$ is 22110; see \cite[A232861]{S} for a list of the first 2000 such numbers $n$.
See also \cite[A236193]{S} for a list of the first 10000 suitable primes $q$ in part (ii) of Conjecture \ref{Conj3.7}.
\end{remark}

\begin{conjecture}\label{Conj3.8} {\rm (i) (2013-12-07)} For every $n=2,3,\ldots$, there is a positive integer $k<n$ with $kp_{n-k}+1$ prime.
Also, for any integer $n>2$, there is a positive integer $k<n$ with $kp_{n-k}-1$ prime.

{\rm (ii) (2013-12-11)} Let $n>5$ be an integer. Then $p_kp_{n-k}-6$ is prime for some $0<k<n$.
\end{conjecture}

\begin{remark}\label{Rem3.8} See \cite[A233296 and A233529]{S} for related data and graphs.
By Chen's work \cite{C}, there are infinitely many primes $p$ with $p+6$ a product of at most two primes.
\end{remark}

\begin{conjecture}\label{Conj3.9} {\rm (2013-11-23)} Let $n>6$ be an integer. Then $p_k+p_{n-k}-1$ is prime for some $k=1,\ldots,n-1$.
Also, $p_k^2+p_{n-k}^2-1$ is prime for some $k=1,\ldots,n-1$.
\end{conjecture}

\begin{remark}\label{Rem3.9} See \cite[A232465]{S} for a related sequence.
\end{remark}

\begin{conjecture}\label{Conj3.10} {\rm (2013-12-10)} {\rm (i)} For any integer $n>3$, there is a positive integer $k<n$ such that $p_k^2+4p_{n-k}^2$ is prime.

{\rm (ii)} Let $n>10$ be an integer. Then there is a positive integer $k<n$ with $p_k^3+2p_{n-k}^3$ prime. Also, $p_k^3+2p_{n-k}^2$ is prime for some $0<k<n$.
\end{conjecture}

\begin{remark}\label{Rem3.10} See \cite[A233439]{S} for a sequence related to part (i). In 2001 Heath-Brown \cite{HB} proved that there are infinitely many primes of the form $x^3+2y^3$ with $x,y\in\Z$.
\end{remark}

\begin{conjecture}\label{Conj3.11} {\rm (2014-03-01)} {\rm (i)} If a positive integer $n$ is not a divisor of $6$, then $p_q^2+(p_n-1)^2$ is prime for some prime $q<n$.
Also, for any positive integer $n\not=1,2,9$, there is a prime $q<n$ with $(p_q-1)^2+p_n^2$ prime.

{\rm (ii)} For every $n=2,3,\ldots$, there is a positive integer $k<n$ with $p_n^3+2p_k^3$ prime.
\end{conjecture}

\begin{remark}\label{Rem3.11} See \cite[A238585]{S} for a sequence related to the first assertion in the conjecture.
\end{remark}

\begin{conjecture}\label{Conj3.12} {\rm (i) (2013-12-05)} Any integer $n>7$ can be written as $k+m$ with $k$ and $m$ positive integers such that
$2^k+p_m$ is prime.

{\rm (ii) (2013-12-06)}  Any integer $n>3$ can be written as $k+m$ with $k$ and $m$ positive integers such that
$k!+p_m$ is prime.
\end{conjecture}

\begin{remark}\label{Rem3.12} See \cite[A233150 and A233206]{S} for related sequences. We have verified parts (i) and (ii) for $n$ up to $3\times10^7$ and $10^7$ respectively.
For example,
for $n=28117716$ we may take $k=81539$ and $m=28036177$ so that $2^k+p_m$ is prime, also $11=4+7$ with $4!+p_7=24+17=41$ prime.
Part (i) was motivated by the author's conjecture (cf. \cite[A231201]{S}) that any integer $n>1$ can be written as a sum of two positive integers $k$ and $m$
with $2^k+m$ prime. We also conjecture that for any integer $n>1$ there is a number $k\in\{1,\ldots,n\}$ with $k!n-1$ (or $k!n+1$) prime.
\end{remark}

\begin{conjecture}\label{Conj3.13} {\rm (i) (2013-12-05)} Any integer $n>2$ can be written as $k+m$ with $m>k>0$ integers such that
$\bi{2k}k+p_m$ is prime.

{\rm (ii) (2014-03-19)} For any integer $n>4$, there exists an integer $1<k<\sqrt n\log n$ such that $p_n+\bi{p_k-1}{(p_k-1)/2}$ is prime.
\end{conjecture}

\begin{remark}\label{Rem3.13} We have verified parts (i) and (ii) for $n$ up to $10^8$ and $10^7$ respectively. See \cite[A233183 and A239451]{S}
for related data and graphs.
\end{remark}

\begin{conjecture}\label{Conj3.14} {\rm (2014-01-21)} For any integer $n\gs20$, there is a positive integer $k<n$ such that
$m=\varphi(k)+\varphi(n-k)/8$ is an integer with $\bi{2m}m+p_m$ prime.
\end{conjecture}

\begin{remark}\label{Rem3.14} This implies that there are infinitely many positive integers $m$ with $\bi{2m}m+p_m$ prime.
(By Stirling's formula, $\bi{2m}m\sim 4^m/\sqrt{m\pi}$ as $m\to+\infty$.)
See \cite[A236241]{S} for the corresponding representation function, and \cite[A236242]{S} for a list of 52 values of $m$ with $\bi{2m}m+p_m$ prime.
For example, when $m=30734$ the number $\bi{2m}m+p_m$ is a prime with 18502 decimal digits.
\end{remark}

\begin{conjecture}\label{Conj3.15} {\rm (i) (2013-12-29)}  Any integer $n>9$ can be written as $k+m$ with $k$ and $m$ positive integers such that
$q=k+p_m$ and $p_q-q+1$ are both prime.

{\rm (ii) (2014-03-05)} Each integer $n>1$ can be written as $k+m$ with $k$ and $m$ positive integers such that
$$p_{p_k}-p_k+1,\ \ p_{p_{2k+1}}-p_{2k+1}+1\ \ \t{and}\ \ p_{p_m}-p_m+1$$
are all prime.

{\rm (iii) (2014-03-06)}  For any positive integer $n$, there is a number $k\in\{1,\ldots,n\}$ such that
$p_{p_k}-p_k+1$ and $p_{p_{kn}}-p_{kn}+1$ are both prime.

{\rm (iv) (2014-03-06)} There are infinitely many primes $q$ with $p_q-q+1$ and $p_{q'}-q'+1$ both prime, where $q'$ is the first prime after $q$.
\end{conjecture}

\begin{remark}\label{Rem3.15} See \cite[A234694, A238766, A238878]{S}
for data and graphs related to parts (i)-(iii). We have verified part (ii) of the conjecture for all $n=2,3,\ldots,10^7$.
See also \cite[A234695]{S} for the first 10000 primes $q$ with $p_q-q+1$ also prime, and \cite[A238814]{S} for the first 10000 primes $q$
with $p_q-q+1$ and $p_{q'}-q'+1$ both prime.
\end{remark}

\begin{conjecture}\label{Conj3.16} {\rm (2014-03-03)}  Let $m>0$ and $n>2m+1$ be integers. If $m=1$ and $2\mid n$, or $m=3$ and $n\not\eq1\ (\mo\ 6)$,
or $m\in\{2,4,5,\ldots\}$, then there is a prime $p<n$ such that $q=\lfloor(n-p)/m\rfloor$ and $p_q-q+1\ (\t{or}\ q^2-2)$ are both prime.
\end{conjecture}

\begin{remark}\label{Rem3.16} In the case $m=1$, this is a refinement of Goldbach's conjecture. When $m=2$, it is stronger
than Lemoine's conjecture which states that any odd number
$n>5$ can be written as $p+2q$ with $p$ and $q$ both prime. Conjecture \ref{Conj3.16} in the case $m>2$ is completely new.
We have verified the conjecture for all $m=1,\ldots,40$ and $n=2m+2,\ldots,10^6$. See \cite[A235189, A238134 and A238701]{S} for related data and graphs.
\end{remark}

\begin{conjecture}\label{Conj3.17} {\rm (2014-01-30)} Any odd number greater than $5$ can be written as a sum of three elements of the set
$$\{q:\ \t{both}\ q \ \t{and}\ p_q-q+1\ \t{are prime}\}.$$
\end{conjecture}

\begin{remark}\label{Rem3.17} This is stronger than the weak Goldbach conjecture finally proved by Helfgott \cite{He}. See \cite[A236832]{S} for the corresponding representation function.
\end{remark}

\begin{conjecture}\label{Conj3.18} {\rm (2014-01-19)} For any integer $n\gs32$, there is a positive integer $k<n-2$ such that
$q=\varphi(k)+\varphi(n-k)/2+1$ and $p_q-q\pm1$ are all prime.
\end{conjecture}

\begin{remark}\label{Rem3.18} See \cite[A236097 and A236119]{S} for related sequences. The conjecture implies that there are infinitely twin prime pairs of the form $p_q-q\pm1$ with $q$ prime.
\end{remark}

\begin{conjecture}\label{Conj3.19} {\rm (2014-01-17)} For any integer $n\gs38$, there is a positive integer $k<n$ such that
$q=\varphi(k)+\varphi(n-k)/3+1$, $r=p_q-q+1$ and $s=p_r-r+1$ are all prime.
\end{conjecture}

\begin{remark}\label{Rem3.19} See \cite[A235924 and A235925]{S} for related sequences. The conjecture implies that there are infinitely primes $q$
with $r=p_q-q+1$ and $s=p_r-r+1$ both prime.
\end{remark}

\begin{conjecture}\label{Conj3.20} {\rm (2014-01-17)} For each $m=2,3,\ldots$, there is a prime chain $q_1<\ldots<q_m$ of length $m$
such that $q_{k+1}=p_{q_k}-q_k+1$ for all $0<k<m$.
\end{conjecture}

\begin{remark}\label{Rem3.20} For such chains of length $m=4,5,6$, see \cite[A235934, A235935 and A235984]{S}. We also have some other conjectures similar to Conjecture \ref{Conj3.20},
see, e.g., \cite[A236066 and A236481]{S}.
\end{remark}

\begin{conjecture}\label{Conj3.21} {\rm (i) (2014-03-06)} For any integer $n>5$, there is a positive integer $k<n$ such that $2k-1$ and $p_{kn}+kn$ are both prime.

{\rm (ii) (2014-01-04)} Any integer $n>8$ can be written as $k(k+1)/2+m$ with $k$ and $m$ positive integers
such that $p_{k(k+1)/2}+\varphi(m)$ is prime.
\end{conjecture}

\begin{remark}\label{Rem3.21} See \cite[A238881 and A235061]{S} for related sequences. Motivated by part (i), we conjecture that
for any integer-valued polynomial $P(x)$ with positive leading coefficient
there are infinitely many positive integers $n$ with $p_n+2P(n)$ prime. When $P(x)$ is constant, this reduces to a conjecture of de Polignac.
\end{remark}

\begin{conjecture}\label{Conj3.22} {\rm (2013-12-16)} {\rm (i)} Any integer $n>100$ can be written as $k^2+m$ with $k$ and $m$ positive integers such that
$\varphi(k^2)+p_m$ is prime.

{\rm (ii)} If an integer $n>6$ is not equal to $18$, then it can be written as $k^2+m$
with $k$ and $m$ positive integers such that $\sigma(k^2) + p_m-1$ is prime, where $\sigma(j)$ is the sum of all positive divisors of $j$.
\end{conjecture}

\begin{remark}\label{Rem3.22} See \cite[A236548]{S} for a sequence related to part (i). Conjecture \ref{Conj3.22} was motivated by the author's conjecture (cf. \cite[A233544]{S}) that
any integer $n>1$ can be written as $k^2+m$ with $\sigma(k^2)+\varphi(m)$  prime, where $k$ and $m$ are positive integers with $m\gs k^2$.
\end{remark}

\begin{conjecture}\label{Conj3.23} {\rm  (2014-02-01)} {\rm (i)} For any integer $n>13$, there is a prime $q<n$ such that $q+2$ and $p_{n-q}+q+1$ are both prime.

{\rm (ii)} If a positive integer $n$ is not a divisor of $12$, then there is a prime $q<n$ such that $3(p_{n-q}+q)-1$ and $3(p_{n-q}+q)+1$
are twin prime.
\end{conjecture}

\begin{remark}\label{Rem3.23}  See \cite[A236831 and A182662]{S} for related sequences.
\end{remark}

\begin{conjecture}\label{Conj3.24} {\rm (2014-05-22)} {\rm (i)} Any integer $n>3$ can be written as $a+b$ with $a$ and $b$ in the set
$$\{k>0:\ \mbox{the inverse of}\ k\ \mbox{mod}\ p_k\ \mbox{is prime}\},$$
where the inverse of $k$ mod $p_k$ refers to the unique $x\in\{1,\ldots,p_k-1\}$ with $kx\eq1\pmod{p_k}$.

{\rm (ii)} Every $n=2,3,4,\ldots$ can be written as $a+b$ with $a$ and $b$ in the set
$$\{k>0:\ k\ \mbox{is a primitive root modulo}\ p_k\}.$$
\end{conjecture}

\begin{remark}\label{Rem3.24} See \cite[A242753 and A242748]{S} for related sequences. We have verified parts (i) and (ii)
for $n$ up to $10^8$ and $3\times10^5$ respectively. We also conjecture that for any prime $p>5$ there is a positive square $k^2<p$ such that the inverse of $k^2$ mod $p$
is prime (cf. \cite[A242425]{S}), and that any integer $n>7$ can be written as $k+m$ with $k,m\in\{2,3,\ldots\}$ such that the least positive residue of $p_k$ modulo $k$
is prime and the least positive residue of $p_m$ modulo $m$ is a square (cf. \cite[A242950]{S}).
\end{remark}

\begin{conjecture}\label{Conj3.25} {\rm (2014-06-01)} Let $n>6$ be an integer. Then there is a prime $p<n$ such that $pn$ is a primitive root modulo $p_n$.
Also, there is a prime $q<n$ such that $q(n-q)$ is a primitive root modulo $p_n$.
\end{conjecture}

\begin{remark}\label{Rem3.25} See \cite[A243164 and A243403]{S} for related data and graphs. We have verified Conjecture \ref{Conj3.25} for all $n=7,\ldots,2\times 10^5$.
\end{remark}

\section{On primes related to partition functions}
\setcounter{lemma}{0}
\setcounter{theorem}{0}
\setcounter{corollary}{0}
\setcounter{remark}{0}
\setcounter{equation}{0}
\setcounter{conjecture}{0}

For $n=1,2,3,\ldots$, let $p(n)$ denote the number of ways to write $n$ as a sum of positive integers with the order of addends ignored.
The function $p(n)$ is called the {\it partition function}.
For each positive integer $n$, let $q(n)$ denote the number of ways to write $n$
as a sum of {\it distinct} positive integers with the order of addends ignored.
The function $q(n)$ is usually called the {\it strict partition function}.
It is known that
$$p(n)\sim\f{e^{\pi\sqrt{2n/3}}}{4\sqrt3 n}\ \ \t{and}\ \ q(n)\sim\f{e^{\pi\sqrt{n/3}}}{4(3n^3)^{1/4}}\ \quad\t{as}\ \ n\to+\infty$$
(cf. \cite{HR} and \cite[p.\,826]{AS}). So both $p(n)$ and $q(n)$ grow eventually faster than any polynomial in $n$.

\begin{conjecture}\label{Conj4.1} {\rm (i) (2014-02-27)} Let $n$ be any positive integer. Then $p(n)+k$ is prime for some $k=1,\ldots,n$. Also,
$q(n)+k$ is prime for some $k=1,\ldots,n$.

{\rm (ii) (2014-02-28)} Let $n>1$ be an integer. Then $p(n)+p(k)-1$ is prime for some $0<k<n$, and $p(k)+q(n)$ is prime for some $0<k<n$.
Also, for any integer $n>7$, there is a positive integer $k<n$ with $n+p(k)$ prime.

{\rm (iii) (2014-03-12)} Let $n>1$ be an integer. Then there exists a number $k\in\{1,\ldots,n-1\}$
such that $kp(n)(p(n)-1)+1$ is prime. Also, we may replace $kp(n)(p(n)-1)+1$ by
$p(k)p(n)(p(n)-1)+1$ or $p(k)p(n)(p(n)+1)-1$.
\end{conjecture}

\begin{remark}\label{Rem4.1} See \cite[A238457, A238509, A239209 and A239214]{S} for related data and graphs. For part (i) or part (ii),
we have verified the first assertion for $n$ up to $1.5\times10^5$. For part (iii), we have verified the first assertion
for $n$ up to $10^5$. Conjecture \ref{Conj4.1} might be helpful in finding large primes.
\end{remark}

\begin{conjecture}\label{Conj4.2} {\rm (2014-02-27)} {\rm (i)} For any integer $n>2$, there is a prime $q<n$ with $2p(n-q)+1$ prime.
Also, for every $n=4,5,\ldots$, there is a prime $q<n$ with $2p(n-q)-1$ prime.

{\rm (ii)} For each integer $n>2$, there is a prime $p<n$ with $q(n-p)+1$ prime. Also,
for any integer $n>6$, there is a prime $p<n$ with $q(n-p)-1$ prime.
\end{conjecture}

\begin{remark}\label{Rem4.2} This is an analogue of Conjecture \ref{Conj2.20}. We have verified the conjecture for $n$ up to $10^5$.
See \cite[A238458 and A238459]{S} for related sequences.
\end{remark}

\begin{conjecture}\label{Conj4.3} {\rm (i) (2013-12-26)} For any integer $n>127$, there is a positive integer $k<n-2$ such that $p(k+\varphi(n-k)/2)$ is prime.

{\rm (ii) (2013-12-28)} For any integer $n>727$, there is a positive integer $k<n-2$ such that $q=\varphi(k)+\varphi(n-k)/2+1$ and $p(q-1)$ are both prime.
\end{conjecture}

\begin{remark}\label{Rem4.3} Clearly, part (ii) implies that there are infinitely many primes $q$ with $p(q-1)$ prime.
We have verified parts (i) and (ii) for $n$ up to $25000$ and $56000$ respectively. See \cite[A234470, A234567 and A234569]{S} for related data and graphs.
\end{remark}

\begin{conjecture}\label{Conj4.4} {\rm (i) (2014-03-13)} For each integer $n>3$, there is a number $k\in\{1,\ldots,n\}$ with $p(n+k)+1$ prime.
Also, for any integer $n>15$, there is a number $k\in\{1,\ldots,n\}$ with $p(n+k)-1$ prime.

{\rm (ii) (2013-12-26)} Any integer $n>5$ can be written as $k+m$ with $k,m\in\{3,4,\ldots\}$ such that $q(\varphi(k)\varphi(m)/4)+1$ is prime.
\end{conjecture}

\begin{remark}\label{Rem4.4} See \cite[A239232 and A234475]{S} for related data and graphs.
The conjecture implies that there are infinitely many primes of the form $p(n)+1\ (\t{or}\ p(n)-1,\ \t{or}\ q(n)+1)$
with $n$ a positive integer.
\end{remark}

\begin{conjecture}\label{Conj4.5} {\rm (2013-12-29)} Any integer $n>7$ can be written as $k+m$ with $k$ and $m$ positive integers such that $p=p_k+\varphi(m)$
and $q(p)-1$ are both prime. Also, any integer $n>7$ not equal to $15$ can be written as $k+m$ with $k$ and $m$ positive integers such that $p=p_k+\varphi(m)$
and $q(p)+1$ are both prime.
\end{conjecture}

\begin{remark}\label{Rem4.5} This implies that there are infinitely many primes $p$ with $q(p)-1$ (or $q(p)+1$) prime. See \cite[A234615 and A234644]{S} for related data and graphs.
We also conjecture that for any integer $n>14$ there exists a prime $p$ with $n<p<2n$ such that $q(p)+1$ is prime.
\end{remark}

\begin{conjecture}\label{Conj4.6} {\rm (2014-01-07)} For any integer $n\gs60$, there is a positive integer $k<n$ such that $m\pm1$ and $q(m)+1$ are all prime, where
$m=\varphi(k)+\varphi(n-k)/4$.
\end{conjecture}

\begin{remark}\label{Rem4.6} This implies that there are infinitely many positive integers $m$ with $m\pm1$ and $q(m)+1$ all prime.
We have verified the conjecture for $n$ up to $10^5$. See \cite[A235343 and A235344]{S} for related data and graphs.
\end{remark}

\begin{conjecture}\label{Conj4.7} {\rm (2014-01-25)} {\rm (i)} For any integer $n\gs128$,
there is a positive integer $k<n$ such that $r=\varphi(k)+\varphi(n-k)/6+1$ and $p(r)+q(r)$ are both prime.

{\rm (ii)} For every $n=18,19,\ldots$, there is a positive integer $k<n$ such that $m=\varphi(k)/2+\varphi(n-k)/8$ is an integer with $p(m)^2+q(m)^2$ prime.
\end{conjecture}

\begin{remark}\label{Rem4.7} Clearly, part (i) implies that there are infinitely many primes of the form $p(r)+q(r)$ with $r$ prime. And part (ii)
implies that there are infinitely many positive integers $m$ with $p(m)^2+q(m)^2$ prime.
We have verified parts (i) and (ii) for $n$ up to $30000$ and $65000$ respectively. See \cite[A236419, A236412 and A236413]{S} for related data and graphs.
\end{remark}

For any positive integer $n$, $\bar q(n)=p(n)-q(n)$ is the number of ways to write $n$
as a sum of unordered positive integers with some part repeated (or even).

\begin{conjecture}\label{Conj4.8} {\rm (2014-01-25)} {\rm (i)} For any integer $n\gs99$, there is a positive integer $k<n$ such that $p=\varphi(k)/2+\varphi(n-k)/12+1$
and $\bar q(p)$ are both prime.

{\rm (ii)} For any integer $n>3$, there is a positive integer $k<n-2$ such that $q(m)^2+\bar q(m)^2$ is prime, where $m=k+\varphi(n-k)/2$.
\end{conjecture}

\begin{remark}\label{Rem4.8} Clearly, part (i) implies that there are infinitely many primes of the form $\bar q(p)$ with $p$ prime. And part (ii)
implies that there are infinitely many positive integers $m$ with $q(m)^2+\bar q(m)^2$ prime.
See \cite[A236417, A236439 and A236440]{S} for related data and graphs.
\end{remark}

\begin{conjecture}\label{Conj4.9} {\rm (i) (2014-03-12)} Let $n>1$ be an integer. Then the number $kp(n)q(n)\bar q(n)-1$ is prime for some $k=1,\ldots,n$.
Also, $2p(k)p(n)q(n)\bar q(n)+1$ is prime for some $k=1,\ldots,n-1$.

{\rm (ii) (2014-01-26)} Any integer $n>2$ can be written as $k+m$ with $k$ and $m$ positive integers such that $q(k)+\bar q(m)$ is prime.

{\rm (iii) (2013-12-08)} For every $n=2,3,\ldots$, there is a positive integer $k<n$ with $2^k-1+q(n-k)$ prime.
\end{conjecture}

\begin{remark}\label{Rem4.9} See \cite[A239207, A236442 and A233390]{S} for related data and graphs.  We have verified the two assertions in part (i) for $n$
up to $83000$ and $50000$ respectively. We have also checked part (iii) for $n$ up to $2\times10^5$; for example, for $n=147650$ we may take $k=17342$ so that $2^k-1+q(n-k)$ is prime.
\end{remark}

\begin{conjecture}\label{Conj4.10} {\rm (2014-04-24)} {\rm (i)} For any prime $p$, there exists a primitive root $g<p$ modulo $p$ which is also a partition number
$($i.e., $g=p(n)$ for some positive integer $n)$.

{\rm (ii)} For any prime $p>3$, there exists a primitive root $g<p$ modulo $p$ which is also a strict partition number
$($i.e., $g=q(n)$ for some positive integer $n)$.
\end{conjecture}

\begin{remark}\label{Rem4.10} We have verified parts (i) and (ii) for all primes below $2\times10^7$ and $5\times10^6$ respectively;
see \cite[A241504 and A241516]{S} for related data and graphs.
We also conjecture that for any prime $p$
there is a primitive root $g<p$ modulo $p$ with $g-1$ a square (cf. \cite[A239957 and A241476]{S}).
\end{remark}

\end{document}